\documentclass[twoside]{article}
\usepackage[mathscr]{eucal}
\usepackage{amsfonts,amssymb,amsmath,graphicx,hyperref}

\oddsidemargin=\evensidemargin
\makeatletter
\long\def\@makefntext#1{
\protect\noindent \hbox to 3.2pt {\hskip-.9pt  
$^{{\eightrm\@thefnmark}}$\hfil}#1\hfill}		

\def\ps@myheadings{\let\@mkboth\@gobbletwo		
\def\@oddhead{\hbox{}
\rightmark\hfil\eightrm\thepage}   
\def\@oddfoot{}\def\@evenhead{\eightrm\thepage\hfil
\leftmark\hbox{}}\def\@evenfoot{}
\def\sectionmark##1{}\def\subsectionmark##1{}}

\def\ps@plain{\let\@mkboth\@gobbletwo
     \def\@oddhead{}\def\@oddfoot{\eightrm\hfil\thepage
     \hfil}\def\@evenhead{}\let\@evenfoot\@oddfoot}

\newcounter{sectionc}\newcounter{subsectionc}\newcounter{subsubsectionc}
\renewcommand{\section}[1] {\vspace{12pt}\addtocounter{sectionc}{1} 
\setcounter{subsectionc}{0}\setcounter{subsubsectionc}{0}\noindent 
	{\tenbf\thesectionc. #1}\par\vspace{5pt}}
\renewcommand{\subsection}[1] {\vspace{12pt}\addtocounter{subsectionc}{1} 
	\setcounter{subsubsectionc}{0}\noindent 
	{\bf\thesectionc.\thesubsectionc. 
	{\kern1pt \bfit #1}}\par\vspace{5pt}}
\renewcommand{\subsubsection}[1] {\vspace{12pt}
	\addtocounter{subsubsectionc}{1}
	\noindent
	{\tenrm\thesectionc.\thesubsectionc.\thesubsubsectionc.	{\kern1pt 
	\it #1}}\par\vspace{5pt}}
\newcommand{\nonumsection}[1] {\vspace{12pt}\noindent{\tenbf #1}
	\par\vspace{5pt}}

\newcounter{appendixc}
\newcounter{subappendixc}[appendixc]
\newcounter{subsubappendixc}[subappendixc]

\renewcommand{\appendix}[1] {\vspace{12pt}	
	\refstepcounter{appendixc}		
	\setcounter{figure}{0}
	\setcounter{table}{0}
	\setcounter{lemma}{0}
	\setcounter{theorem}{0}
	\setcounter{corollary}{0}
	\setcounter{definition}{0}
	\setcounter{equation}{0}
	\renewcommand{\thefigure}{\Alph{appendixc}.\arabic{figure}}
	\renewcommand{\thetable}{\Alph{appendixc}.\arabic{table}}
	\renewcommand{\theappendixc}{\Alph{appendixc}}
	\renewcommand{\thelemma}{\Alph{appendixc}.\arabic{lemma}}
	\renewcommand{\thetheorem}{\Alph{appendixc}.\arabic{theorem}}
	\renewcommand{\thedefinition}{\Alph{appendixc}.\arabic{definition}}
	\renewcommand{\thecorollary}{\Alph{appendixc}.\arabic{corollary}}
	\renewcommand{\theequation}{\Alph{appendixc}.\arabic{equation}}
	\noindent{\tenbf Appendix \theappendixc #1}\par\vspace{5pt}}

\newcommand{\textlineskip}{\baselineskip=13pt}
\newcommand{\smalllineskip}{\baselineskip=10pt}

\def\abstracts#1#2#3#4{{
	\centering{\begin{minipage}{4.5in}\footnotesize\baselineskip=10pt
	\centerline{ABSTRACT} 
	\parindent=15pt #1\par 
	\parindent=15pt #2\par
	\parindent=15pt #3\par
	\parindent=15pt #4\par
	\end{minipage}}\par}} 

\def\keywords#1{{ 
	\centering{\begin{minipage}{4.5in}\footnotesize\baselineskip=10pt
	{\footnotesize\it Keywords}\/: #1
	\end{minipage}}\par}}

\renewenvironment{thebibliography}[1]
	{\frenchspacing
	 \ninerm\baselineskip=11pt
	 \begin{list}{[\arabic{enumi}]}
	{\usecounter{enumi}\setlength{\parsep}{0pt}
	 \setlength{\leftmargin 19pt}{\rightmargin 0pt}   
	 \setlength{\itemsep}{0pt} \settowidth
	{\labelwidth}{[#1]}\sloppy}}{\end{list}}

\newcounter{itemlistc}
\newcounter{romanlistc}
\newcounter{alphlistc}
\newcounter{arabiclistc}

\newcommand{\fcaption}[1]{
        \refstepcounter{figure}
        \setbox\@tempboxa = \hbox{\footnotesize Fig.~\thefigure. #1}
        \ifdim \wd\@tempboxa > 5in
           {\begin{center}
        \parbox{5in}{\footnotesize\smalllineskip Fig.~\thefigure. #1}
            \end{center}}
        \else
             {\begin{center}
             {\footnotesize Fig.~\thefigure. #1}
              \end{center}}
        \fi}

\newcommand{\tcaption}[1]{
        \refstepcounter{table}
        \setbox\@tempboxa = \hbox{\footnotesize Table~\thetable. #1}
        \ifdim \wd\@tempboxa > 5in
           {\begin{center}
        \parbox{5in}{\footnotesize\smalllineskip Table~\thetable. #1}
            \end{center}}
        \else
             {\begin{center}
             {\footnotesize Table~\thetable. #1}
              \end{center}}
        \fi}

\def\pmb#1{\setbox0=\hbox{#1}
	\kern-.025em\copy0\kern-\wd0
	\kern.05em\copy0\kern-\wd0
	\kern-.025em\raise.0433em\box0}

\def\fnt#1#2{\footnotetext{\kern-.3em
	{$^{\mbox{\scriptsize #1}}$}{#2}}}

\def\fpage#1{\begingroup
\voffset=.3in
\thispagestyle{empty}\begin{table}[b]\centerline{\footnotesize #1}
	\end{table}\endgroup}

\def\runninghead#1#2{\pagestyle{myheadings}
\markboth{{\protect\footnotesize\it{\quad #1}}\hfill}
{\hfill{\protect\footnotesize\it{#2\quad}}}}

\font\tenrm=cmr10
 
\font\tenbf=cmbx10
\font\bfit=cmbxti10 at 10pt
\font\ninerm=cmr9
\font\nineit=cmti9
\font\ninebf=cmbx9
\font\eightrm=cmr8

\newtheorem{theorem}{Theorem}   

\newtheorem{lemma}{Lemma}

\newtheorem{corollary}{Corollary}
\newtheorem{proposition}{Proposition}
\newtheorem{remark}{Remark}
\newtheorem{example}{Example}
\newtheorem{construction}{Construction}

\def\@begintheorem#1#2{\trivlist	
	\item[\hskip\labelsep{\bf #1\ #2.}]} 
\def\@opargbegintheorem#1#2#3{\trivlist
	\item[\hskip\labelsep{\bf #1\ #2\ (#3).}]}

\newenvironment{proof}{\begin{trivlist}
	\item[\noindent]{\it Proof.}}{\quad $\square$\end{trivlist}}

\newenvironment{romanlist2}[1]			
	{\setcounter{romanlistc}{0}		
	 \begin{list}{$($\roman{romanlistc}$)$}	
	{\usecounter{romanlistc}		
	 \leftmargin18pt 
	 \setlength{\parsep}{0pt}
	 \setlength{\itemsep}{0pt}	
	 \settowidth{\labelwidth}{#1}                          
	}}{\end{list}}

\def\qed{\hbox{${\vcenter{\vbox{			
   \hrule height 0.4pt\hbox{\vrule width 0.4pt height 6pt
   \kern5pt\vrule width 0.4pt}\hrule height 0.4pt}}}$}}

\def\theequation{\thesectionc.\arabic{equation}}  

\makeatother

\def\Ps{\mathscr P}
\def\Rs{\mathscr R}
\def\Ss{\mathscr S}

\def\C{\mathbb{C}}

\def\R{\mathbb{R}}
\def\Z{\mathbb{Z}}

\def\card{\operatorname{card}}

\def\link{\operatorname{link}}

\def\ord{\operatorname{ord}}
\def\gr{\operatorname{gr}}

\def\TB{\operatorname{TB}}

\def\Bd{\partial}
\def\eps{\varepsilon}
\def\phi{\varphi}
\def\Ncirc{\hbox{$N$\kern-0.75em\raise0.5ex\hbox{\text{\char'27}}}}
\def\s{\sigma}
\def\sub{\subset}
\def\emptyset{\varnothing}
\def\overstrike#1#2{{\setbox0\hbox{$#2$}\hbox to \wd0{\hss
                         $#1$\hss}\kern-\wd0\box0}}
\def\twobars#1#2#3#4#5#6{\vcenter{\hrule height.#1pt width#2pt
                               \vskip#3pt
                               \hrule height.#4pt width#5pt
                               \vskip#6pt}}
     \def\stroke#1#2#3{\vrule height#1pt width.#2pt depth#3pt}
     \def\connsum{\raise.25ex\hbox{\overstrike\parallel=}}
     \def\bdconnsum{\hskip2pt
                   \stroke83{-1.15}\twobars433431\stroke{4.85}31
                   \hskip2pt}
\def\Connsum{\hskip1pt
             \twobars523523\stroke952\twobars533533\stroke952\twobars523523
             \hskip.5pt}
\def\Bdconnsum{\hskip1pt
             \stroke95{-2.05}\twobars533533\stroke65{1.5}
             \hskip.5pt}

\begin{document}
\setlength{\textheight}{7.7truein}  

\runninghead{Lee Rudolph}{Quasipositive Annuli}

\normalsize\textlineskip
\thispagestyle{empty}
\setcounter{page}{1}

{\vspace*{-2.5cm}\smalllineskip{\flushleft
{\footnotesize Hyper{\TeX}ed arXival version prepared December 2001; 
       originally published in}\\
{\footnotesize Journal of Knot Theory and Its Ramifications 
        Vol.~1, No.~4 (1992) 451--466.}\\
         }}

\vspace*{0.88truein}

\fpage{1}
\centerline{\bf QUASIPOSITIVE ANNULI (CONSTRUCTIONS OF}
\baselineskip=13pt
\centerline{\bf QUASIPOSITIVE KNOTS AND LINKS, IV)}
\vspace*{0.37truein}
\centerline{\footnotesize LEE RUDOLPH}
\baselineskip=12pt
\centerline{\footnotesize\it Department of Mathematics 
and Computer Science}
\baselineskip=10pt
\centerline{\footnotesize\it Clark University, Worcester, 
Massachusetts 01610}

\vspace*{0.21truein} 
\abstracts{The {\it modulus of quasipositivity} $q(K)$ of a 
knot $K$ was introduced as a tool in the knot theory of complex 
plane curves, and can be applied to Legendrian knot theory in symplectic 
topology. It has also, however, a straightforward characterization 
in ordinary knot theory: $q(K)$ is the supremum of the integers 
$f$ such that the framed knot $(K,f)$ embeds non-trivially 
on a fiber surface of a positive torus link. Geometric constructions 
show that $-\infty <q(K)$, calculations with link polynomials that 
$q(K)<\infty $.  The present paper aims to provide sharper lower 
bounds (by optimizing the geometry with {\it positive plats}) and 
more readily calculated upper bounds (by modifying known link 
polynomials), and so to compute $q(K)$ for various classes of knots, 
such as positive closed braids (for which $q(K)=\mu(K)-1)$ and 
most positive pretzels.  As an aside, it is noted that a recent
result of Kronheimer~\&~Mrowka implies that $q(K)<0$ if $K$ is slice.}{}{}{}

\vspace*{10pt}

\vspace*{1pt}\textlineskip	
\section{Review of Background Material}	
\vspace*{-0.5pt}
By default, manifolds are piecewise-smooth, compact, unbounded, 
and oriented; $-X$ denotes $X$ with orientation reversed.

\subsection{Annular Seifert Surfaces and Framed Links}

A surface is {\it annular} if each component is an annulus. 
A {\it Seifert surface} $S\sub S^3$ is a 2-submanifold-with-boundary 
such that each component of $S$ has non-empty boundary. Let 
$K\sub S^3$ be a knot, $f\in \Z$; by $A(K,f)$ denote any annulus 
in $S^3$ with $K\sub\Bd A(K,f)$,  $\link(K,\Bd A(K,f)\setminus K)=-f$.  
Up to ambient isotopy, $A(K,f)$ depends only 
on $K$ and $f$,  and $A(K,f)=A(-K,f)=-A(K,f)$.  Let $L\sub S^3$ 
be a link with components $K_i$,  $f:L\to\Z$ a {\it framing} of $L$ 
(i.e., a continuous function); by $A(L,f)$ denote 
any union of pairwise disjoint annuli $A(K_i,f(K_i))$.  Any annular 
Seifert surface has the form $A(L,f)$.  A framed link ($L,f)$ is 
{\it embedded} on a Seifert surface $S$ if $L\sub S$ 
and the regular neighborhood $N_S(L)$ is $A(L,f)$. 
 
\begin{lemma}
\label{sum lemma}
If $(K_k,f_k)$ is embedded on $S_k$, $k=1,2$,
then the connected sum $(K_1\connsum K_2,f_1 + f_2)$ embeds on a
boundary-connected sum $S_1 \bdconnsum \,\eps S_2\;(\eps\in \{ +,- \})$.
\end{lemma}

\begin{figure}[htbp]
\vspace*{13pt}
\setlength{\unitlength}{1in}
\centering
\begin{picture}(1.375,1)(0,0)
\includegraphics[width=1.375in]{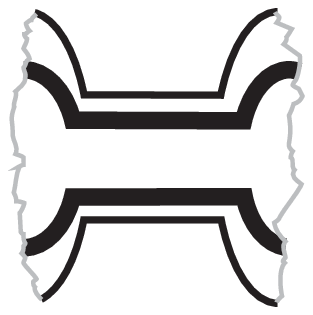}
\put(-1.6,1.0){$K_1$}\put(0.05,1.0){$K_2$}
\put(-1.2,0.65){$S_1$}\put(-0.35,0.65){$\pm S_2$}
\put(-1.6,0.2){$K_1$}\put(0.05,0.2){$K_2$}
\end{picture}
\vspace*{13pt}
\fcaption{\label{4FIG1}}
\end{figure}

\begin{proof} Fig.~\ref{4FIG1}.
\end{proof}

\begin{remark} If $S$ is a Seifert surface with connected 
boundary $K$, then $(K,0)$ and $(-K,0)$ both embed on $S$.  If $K$ 
is non-invertible, that is, not isotopic to $-K$,  then $(K\connsum{-K},0)$ 
embeds on $S \bdconnsum{-S}$ but need not embed on $S \bdconnsum S$.
\end{remark}

\subsection{Closed Braids, Plats, Band Representations, 
and Braided Surfaces}

\begin{figure}[htbp]
\vspace*{13pt}
\setlength{\unitlength}{1in}
\centering
\begin{picture}(1.25,1.25)(0,0)
\includegraphics[width=1.375in]{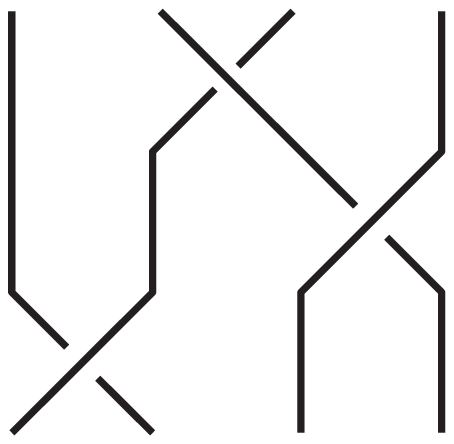}
\end{picture}
\vspace*{13pt}
\fcaption{The braidword diagram of $\mathbf b = 
(\s_1,\s_3,\s_2^{-1})$ ($n = 4$).\label{4FIG2}}
\end{figure}

For present purposes, we may define the {\it $n$-string braid group} 
by its {\it standard presentation} 
$$
B_n:=
\text{gp}\left( \s _i,1\leq i\leq n-1\left| 
\genfrac{}{}{0pt}{}{[\s_i, \s_j]=\s_j^{-1}\s_i,}{[\s_i, \s_j]=1,}
\genfrac{}{}{0pt}{}{|i-j|=1}{|i-j|\not=1}
\right.\right). 
$$
A {\it braidword} in $B_n$ is a $k$-tuple 
$\mathbf b=(\s_{i(1)}^{\eps(1)},\dots,\s_{i(k)}^{\eps(k)})$,  
$\eps(s)\in \{1,-1\}$; the {\it braid} of $\mathbf b$ is 
$\beta(\mathbf b):=\s_{i(1)}^{\eps(1)}\dotsm\s_{i(k)}^{\eps(k)}$.  
It is usual to picture a braidword by a {\it braidword diagram} with 
$2n$ loose ends, $n$ at the top and $n$ at the bottom, 
cf.~Fig.~\ref{4FIG2}.

\begin{figure}[htbp]
\vspace*{13pt}
\setlength{\unitlength}{1in}
\centering
\begin{picture}(3.5,3.0)(0,0)
\includegraphics[width=3.5in]{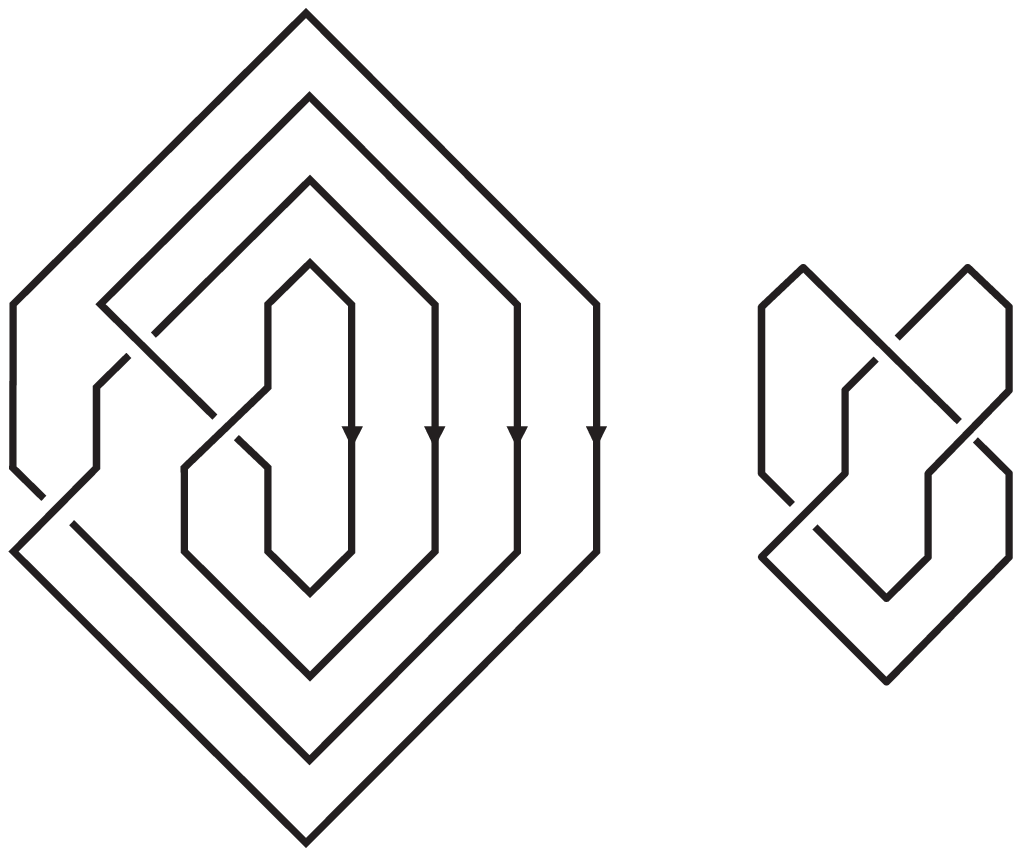}
\put(-0.85,2.15){$\pi_{\sqcap}=(12)(34)$}
\put(-0.85,0.4){$\pi_{\sqcup}=(14)(23)$}
\end{picture}
\vspace*{13pt}
\fcaption{The closed braid $\widehat{\beta}$ and a  
$\Pi$-plat $\beta ^\Pi $ for $\beta=\s_1\s_3\s_2^{-1}\in B_4$.
\label{4FIG3}}
\end{figure}

Let $\beta\in B_n$.  Fig.~\ref{4FIG3}(1) 
indicates the construction of a link 
$\widehat{\beta} \sub S^3$, called the {\it closed braid} 
of $\beta$: let $\mathbf b$ be any braidword with $\beta =\beta(\mathbf b)$; 
then a link diagram for 
$\widehat{\beta} (\mathbf b):=\widehat{\beta(\mathbf b)}$ 
consists of the braidword diagram of $\mathbf b$, together with $n$ arcs 
joining the loose ends at the bottom to those at the top 
so as to create no new crossings. The closed braid 
$\widehat{\beta}$ depends only on $\beta $ (not on $\mathbf b)$,  
is well-defined up to ambient isotopy, and has the canonical orientation 
indicated in Fig.~\ref{4FIG3}(1).

Let $n$ be even. A {\it plat-plan} is a pair 
$\Pi=(\pi_\sqcup,\pi_\sqcap)$ of permutations in ${\Ss}_n$
each of which is the product of $n/2$ disjoint transpositions 
such that 
$s<t<\pi(s)\Rightarrow s< \pi(t)< \pi(s)$.  
Fig. 3(2) indicates how, given a plat-plan $\Pi$ and $\beta\in B_n$,
to construct an unoriented link $\beta^\Pi\sub S^3$, 
called the {\it $\Pi$-plat} of $\beta $: 
let $\mathbf b$ be any braidword with 
$\beta=\beta(\mathbf b);$ then a link diagram for 
$\beta^\Pi(\mathbf b):=\beta(\mathbf b)^\Pi$ 
consists of the braidword diagram of $\mathbf b$,  together 
with $n/2$ arcs joining the loose ends at the bottom 
according to $\pi_\sqcup $, and $n/2$ arcs joining 
the loose ends at the top according to $\pi_\sqcap $,  all so as 
to create no new crossings. The $\Pi$-plat $\beta^\Pi$ 
depends only on $\beta$ (not on $\mathbf b)$ and is well-defined up 
to ambient isotopy; it has no canonical orientation. Where possible, 
the particular plat-plan $\Pi$ is suppressed and a $\Pi$-plat
is simply called a {\it plat}. 
 
\begin{remark} The {\it writhe} (i.e., algebraic 
crossing number) of a knot diagram (as opposed to a multi-component 
link diagram) derived as above from a braidword diagram of $\mathbf b$,  
whether for the closed braid or the $\Pi$-plat of $\beta(\mathbf b)$,  
depends only on $\beta(\mathbf b)$; in the case of the closed braid, 
this writhe equals $e(\beta(\mathbf b))$, 
the exponent sum of $\beta (\mathbf b)$ (with respect to the
standard generators $\s_i$ of $B_n$), but there seems to be 
no similarly neat expression for the writhe of a plat.
\end{remark} 

\begin{figure}[htbp]
\vspace*{13pt}
\setlength{\unitlength}{1in}
\centering
\begin{picture}(2.75,2.5)(0,0)
\includegraphics[width=2.5in]{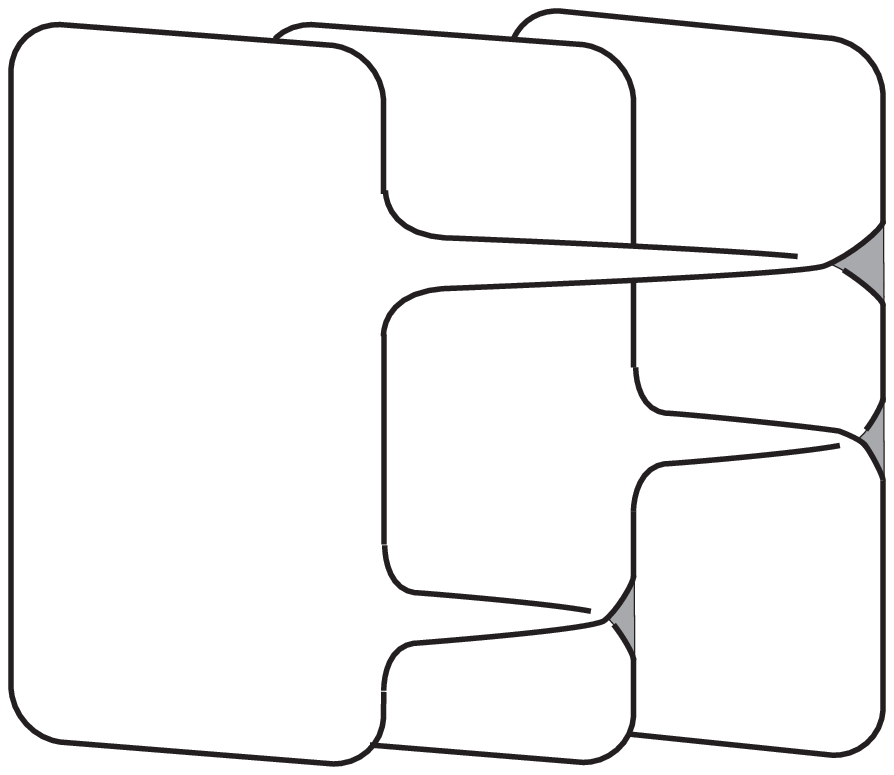}
\end{picture}
\vspace*{13pt}
\fcaption{A braided Seifert surface $S(\mathbf b)$,  
$\mathbf b=(\s_{1,2}^{\phantom{1}},\s_{2,3}^{-1},\s_{1,3}^{\phantom{1}})$.\
\label{4FIG4}}
\end{figure}

A {\it positive embedded band} in $B_n$ is one of the 
$\binom{n}{2}$ braids 
$$
\s_{i,j}:=
(\s_i\dotsm\s_{j-2})\s_{j-1}(\s_i\dotsm\s_{j-2})^{-1},
\quad1\leq i<j\leq n;
$$ 
a {\it negative embedded band} is the inverse of a positive embedded band.
An {\it embedded band representation} 
of {\it length} $k$ in $B_n$ is a $k$-tuple $\mathbf b=:
(b(1),\dotsm,b(k))$ of embedded bands; as with braidwords 
(which are embedded band representations, 
since $\s_i=\s_{i,i+1}$),  write 
$\beta(\mathbf b):=b(1)\dotsm b(k)$,  
$\widehat{\beta}(\mathbf b):=\widehat{\beta(\mathbf b)}$,  
$\beta^\Pi(\mathbf b):=\beta(\mathbf b)^\Pi$.  
There is a straightforward construction\cite{1,2,3,4}
of a {\it braided Seifert surface}
$$
S(\mathbf b)=\bigcup^n_{s=1}h_i^{(0)}\cup\bigcup^k_{t=1}h_j^{(1)}
$$
given as the union of $n$ $0$-handles and $k$ $1$-handles, 
with $\Bd S(\mathbf b)=\widehat\beta(\mathbf b)$ 
(cf.~Fig.~\ref{4FIG4}). 

Call $\mathbf b$ {\it weakly annular} if $S(\mathbf b)$ is annular, 
and {\it annular} if also every $0$-handle $h_s^{(0)}$ of $S(\mathbf b)$ 
is attached to two $1$-handles $h_{i(s)}^{(1)}, h_{j(s)}^{(1)}$ 
(rather than to only a single $1$-handle). 

\subsection{Quasipositivity}

An embedded band representation $\mathbf b$ is {\it quasipositive} 
if each $b(t)$ is positive. A Seifert surface $S$ is {\it quasipositive} 
if, for some quasipositive $\mathbf b$,  $S$ is ambient isotopic 
to $S(\mathbf b)$.  In this case, $\mathbf b$ may always be taken to be 
such that no $0$-handle of $S(\mathbf b)$ is attached to precisely 
one $1$-handle; in particular, if $S$ is annular, then $\mathbf b$ 
may also be taken to be annular (rather than merely weakly annular).
 
\begin{example}
\label{positive closed braids}
A braidword $\mathbf p$ is quasipositive if 
and only if it is {\it positive} (i.e., $\eps(t)=1$ 
for all $t$).  A positive closed braid $\widehat\beta(\mathbf p)$ 
has many special properties (\cite{5,6,7}), 
among them that each 
component of $S(\mathbf p)$ is a fiber surface, so the split components 
of $(S^3,\widehat\beta(\mathbf p))$ are fibered links. 
In particular, for $m,n>0$,  the torus link $O\{m,n\}$ 
is the non-split positive closed braid $\widehat\beta(\mathbf o\{m,n\})$,  
where $o\{m,n\}(t)$ $:=$ $\s_i\in B_m$ for $t\equiv i\pmod{m-1}$,  
$1\leq t\leq n(m-1)$; $O\{ m,n\}$ is the link of the singularity 
at the origin of the complex plane curve $\{(z,w)\in \C^2:z^m+w^n=0\}$, 
as well as the link at infinity of the same curve, so its fiber surface 
$S(\mathbf o\{m,n\})$ can be realized as 
$\{(z,w)\in S^3:z^m+w^n\geq 0\}\sub S^3
:=\{(z,w)\in\C^2:\left|z\right|^2+\left|w\right|^2=1\}$. 
\end{example}

Let $S$ be a surface. A {\it full} subset $X$ of $S$ is 
one such that no component of $S\backslash X$ is contractible (e.g., 
a collar of $\Bd S$ is full iff no component of $S$ 
is $D^2$; a simple closed curve on $S$ is full iff it
does not bound a disk on $S$). 
 
\begin{theorem} (\cite{8}) If $S$ is a Seifert surface, then 
the following are equivalent. 
\begin{romanlist2}{(iv)}
\item $S$ is quasipositive. 
\item $S$ is a full subsurface of a fiber surface 
of a positive closed braid. 
\item $S$ is a full 
subsurface of a fiber surface of a positive torus link. 
\item $S$ is a full subsurface of a fiber surface of 
$O\{ n,n\}$ for some $n>0$.  
\end{romanlist2}
\end{theorem}
 
\begin{corollary}
\label{full subsurface of qp is qp} 
A full subsurface of a quasipositive 
Seifert surface is quasipositive. 
\end{corollary}

The {\it modulus of quasipositivity} of a knot $K$ is 
\begin{align*}\label{definition of q}
q(K)&:=\sup \{f\in \Z : A(K,f)\text{ is quasipositive}\}\\
    &\phantom{:}=\sup\{f:A(K,f)\text{ is full on 
                                     some quasipositive surface}\}\\
    &\phantom{:}=\sup\{f: \text{ for some $n$, $A(K,f)$ is full on 
                                     $S(\mathbf o\{n,n\})$}\}
\end{align*}
(the equalities following from 
Corollary~\ref{full subsurface of qp is qp}).  For any $K$, $-\infty<q(K)$ 
by \cite{3} and $q(K)<\infty $ by \cite{9}.

\begin{proposition}
\label{quasipositive annuli} 
(\cite{3}) For any knot $K$, if $f\leq q(K)$ 
then $A(K,f)$ is quasipositive.
\end{proposition}
 
\begin{proof} For $n\geq 1$, a collar of a component of the boundary of 
a fiber surface of $O\{ n+1,n+1\} $ is an annulus 
$A(O,-n)$; the proposition follows from Lemma~\ref{sum lemma}, 
Corollary~\ref{full subsurface of qp is qp}, and the 
observation that any boundary-connected sum of quasipositive Seifert surfaces 
is quasipositive (e.g., 
$S(\s_1,\s_1,\s_2,\s_2,\s_2)\bdconnsum S(\s_1,\s_1,\s_1)$ is either
$S(\s_1,\s_1,\s_2,\s_2,\s_2,\s_3,\s_4,\s_4,\s_4)$ or
$S(\s_1,\s_1,\s_2,\s_2,\s_2,\s_{1,4},\s_4,\s_4,\s_4)$). 
\end{proof}  
 
\begin{corollary}\label{subadditivity of q} For any knots $K_k$,  
$q(\Connsum_k K_k)\geq \sum _k q(K_k)$. 
\end{corollary}
 
\begin{proof} By Lemma~\ref{sum lemma}, $(\Connsum_k K_k ,\sum _k f_k)$ 
embeds on $\Bdconnsum_k A(K_k,f_k)$ 
(since $A(K,f)=-A(K,f)$ for all $(K,f)$).  
\end{proof} 
 
In Proposition~\ref{subadditivity of q+1} this is improved to 
$q(\Connsum_k K_k)+1\geq \sum _k(q(K_k)+1)$.   

\subsection{Link Polynomials}

\begin{figure}[htbp]
\setlength{\unitlength}{1in}
\centering
\begin{picture}(3,0.85)(-0.01,-0.10)
\includegraphics[width=3in]{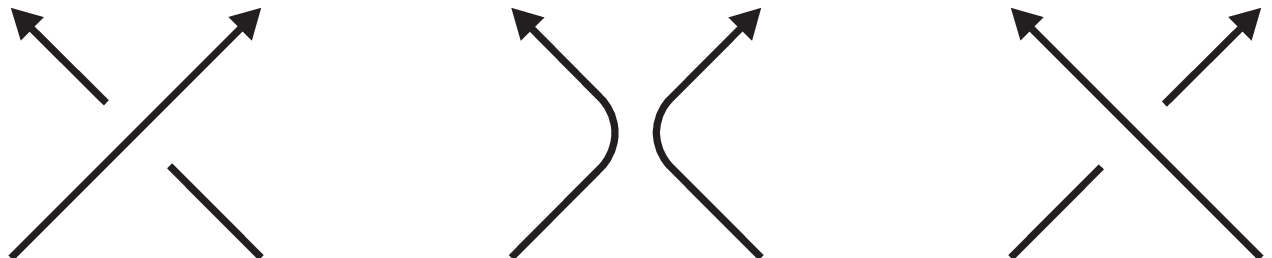}
\end{picture}
\put(-2.75,0){$L_{+}$}\put(-1.6,0){$\thinspace L_0$}
\put(-.45,0){$\quad\negthinspace L_{-}$}
\fcaption{\label{4FIG5}}
\end{figure}

\begin{figure}[htbp]
\setlength{\unitlength}{1in}
\centering
\begin{picture}(1.8,.85)(0.0,.10)
\includegraphics[width=1.8in]{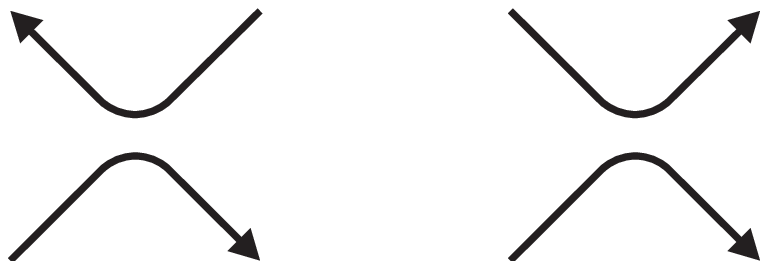}
\put(-1.6,-.10){\mbox{\scriptsize{$(1)$}}}
\put(-0.4,-.10){\mbox{\scriptsize{$(2)$}}}
\end{picture}
\fcaption{\label{4FIG6}}
\end{figure}

If $L_{+}$, $L_0$, and $L_{-}$ are links with diagrams which are identical 
except as indicated in Fig.~\ref{4FIG5}, and the visible crossing in $L_{+}$ 
involves one component (resp., two components), then we say we are 
in case $1$ (resp., case $2$), we let $p$ (resp., $q$) be the linking number 
of the right-hand visible component of $L_0$ with the rest of $L_0$ 
(resp., the linking number of the lower visible component of $L_{+}$ 
with the rest of $L_{+}$), and we denote by $L_\infty$ the link 
indicated in Fig.~\ref{4FIG6}(1) (resp., Fig.~\ref{4FIG6}(2)). 

The {\it oriented} (or {\it FLYPMOTH} \cite{10,11}) 
and {\it semi-oriented} (or {\it Kauffman} \cite{12}) polynomials 
$P_L(v,z)\in\Z[v^{\pm1},z^{\pm1}]$ and 
$F_L(a,x)\in\Z[a^{\pm1},x^{\pm1}]$ 
of a link $L$ can 
be defined recursively as follows: 
\begin{gather*}
P_O(v,z)=1=F_O(a,x) \text{ if $L=O$ is an unknot,}\\
P_{L_{+}}(v,z)=vzP_{L_0}(v,z)+v^2P_{L_{-}}(v,z),\\
F_{L_{+}}(a,x)=a^{-1}xF_{L_0}(a,x)-a^{-2}F_{L_{-}}(a,x)+
\begin{cases}
a^{-4p-1}x F_{L_\infty}(a,x)& \text{ in case $1$,}\\
a^{-4q+1}x F_{L_\infty}(a,x)& \text{ in case $2$}\\
\end{cases}
\end{gather*}
(the use of the names ``oriented'' and ``semi-oriented'' follows 
Lickorish \cite{13}; the use of the variables $v$ and $z$ follows Morton 
\cite{14}).

For any ring $\Rs$,  for any Laurent polynomial $H(s)\in \Rs[s^{\pm 1}]$, 
write $\ord_s H(s):=\sup\{ n\in\Z:s^{-n}H(s)\in\Rs[s]\sub\Rs[s^{\pm 1}]\}$, 
$\deg_s H(s):=-\ord_s H(s^{-1})$.  
Let $c(L)$ be the number of components of $L$.  Easy inductions establish
the following estimates.
 
\begin{lemma}
\label{first estimate} For every link $L$, $\ord_z P_L\geq 1-c(L)$ and 
$\ord_x F_L\geq 1-c(L)$.
\end{lemma} 
 
\begin{corollary}\label{well-defined}  
$\left.(z^{c(L)-1}P_L(v,z))\right|_{z=0}$
and $\left.(x^{c(L)-1}F_L(a,x)\right|_{x=0}$
are well-defined link invariants (in $\Z[v^{\pm 1}]$ and 
$\Z[a^{\pm 1}]$, respectively). 
\end{corollary}
 
Let $R_L(v):=\left.(z^{c(L)-1}P_L(v,z))\right|_{z=0}$. 
 
\begin{lemma} (cf.~\cite{15}) $\left.((-\sqrt{-1}x)^{c(L)-1}F_L(a,x)\right|_{x=0}=
R_L(\sqrt{-1}a^{-1})$.
\end{lemma}

\begin{lemma}
\label{R-polynomial}
$R_L(v)$ can be calculated recursively as follows:
\begin{gather*}
R_O(v)=1,\\
R_{L_{+}}(v)=(2-k)vR_{L_0}(v)+v^2R_{L_{-}}(v)\text{ in case $k$ ($k=1,2$).}
\end{gather*}
\end{lemma}
 
\begin{corollary}
\label{R of disjoint sum}
If links $L_1$ and $L_2$ are disjoint (i.e., 
if $L_1\cup L_2$ is a link), then  
$$
R_{L_1\cup L_2}(v)=(v^{-1}-v)v^{2\link(L_1,L_2)}R_{L_1}(v)R_{L_2}(v).
$$
\end{corollary}
 
If $L$ has components $K_i$,  then its {\it total linking}
is $\tau(L):=\sum_{i<j}\link(K_i,K_j)$.  If $f$ is a 
framing of $L$,  then the {\it total framing} of the framed 
link $(L,f)$ is $\phi(L,f):=\sum_i f(K_i)$.  Define the 
{\it framed polynomial} to be
$$
\{L,f\}(v,z):=(-1)^{c(L)}(1
+(v^{-1}-v)z^{-1} \sum_{L'}(-1)^{c({L'})}P_{\Bd A(L',f|{L'})})
$$
where $L'$ runs through the non-empty sublinks of $L$.  

\begin{proposition} $\{L,f\}=v^{-2\phi(L,f)}\{L,0\}$. 
\end{proposition}

The framed polynomial provides a bridge between the oriented and semi-oriented 
polynomials, as the following result (proved in \cite{9}) makes plain.
 
\begin{theorem}
\label{congruence theorem}
$(1+(v^{-2}+v^2)z^{-2})F_L(v^{-2},z^2)
\equiv v^{4\tau(L)}\{L,0\}(v,z)\pmod{2}$.  
\end{theorem}

Let $F^*_L:=(F_L\pmod{2})\in(\Z/2\Z)[a^{\pm1},x^{\pm 1}]$,
$G^k_L(a):=\left.(x^{1-c(L)}F^*_L(a,x))\right|_{x=k}
\in(\Z/2\Z)[a^{\pm 1}]$ for $k=0,1$.  Thus 
$G^0_L(a)=R_L(a^{-1}) \pmod{2}$, and can be 
calculated using the formulas in Lemma~\ref{R-polynomial} reduced mod $2$.

\begin{lemma}
\label{G-polynomial}
The polynomial $G^1_L$ can be calculated recursively as follows:
\begin{gather*}
G^1_{O}(a)=1,\\
G^1_{L{_{+}}}(a)=a^{-2}G^1_{L{_{-}}}(a)+a^{-1}G^1_{L{_0}}(a)+
\begin{cases}
{a^{-4p-1}G^1_{L{_\infty}}(a) \text{ in case $1$,}}\\
{a^{-4q+1}G^1_{{L_\infty}}(a) \text{ in case $2$.}} 
\end{cases}
\end{gather*}
\end{lemma}

\section{Lower Bounds for the Modulus of Quasipositivity}

\subsection{Summary of Results}

In this section we use {\it fences} to show that, 
in contrast to positive closed braids (some special properties of which
were mentioned above in Example~\ref{positive closed braids}), 
positive plats are not at all special: every (unoriented) link 
is ambient isotopic to a plat of a positive braidword 
(Construction~\ref{fence to positive plat}). 
We then express $q(K)$ in terms of positive plat 
realizations of $K$ (Corollary~\ref{q from positive plats}), 
thus bounding $q(K)$ from below (Theorem~\ref{bound q from below}).

\subsection{Fences }

\begin{figure}[htbp]
\vspace*{13pt}
\setlength{\unitlength}{1in}
\centering
\begin{picture}(3.5,2.4)(0,0)
\includegraphics[width=3.5in]{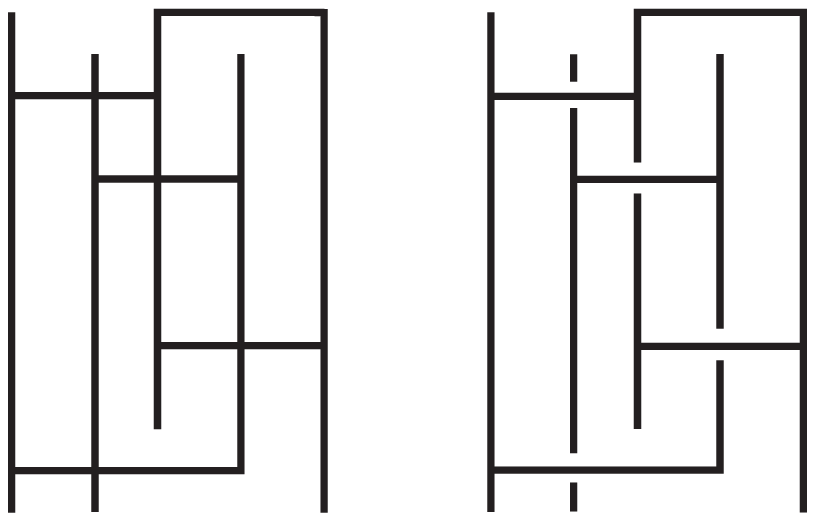}
\end{picture}
\vspace*{13pt}
\fcaption{A fence $\Phi$ and a graph $\gr(\Phi)$.\label{4FIG7}}
\end{figure}

A {\it post} is a vertical segment $\{ x\}\times[a,b]\sub\R^2$,  
$a<b$.  A {\it wire} is a horizontal segment $[c,d]\times\{y\}\sub\R^2$,
$c<d$.  A {\it fence} $\Phi\sub\R^2$ is the union of $n\geq 1$ posts 
with pairwise distinct abscissae and $k\geq 0$ wires 
with pairwise distinct ordinates, 
such that both endpoints of each wire lie on 
posts; cf.~Fig.~\ref{4FIG7}(1). 

Let $X(\Phi)$ denote the set of abscissae 
of posts of $\Phi $, $Y(\Phi)$ the set of ordinates of wires of 
$\Phi$; let $x_1<\dots<x_n$ be the elements of $X(\Phi)$,  
$y_1<\dots<y_k$ the elements of $Y(\Phi)$.  
For $1\leq t\leq k$, define $i(t)$ and $j(t)$ by the requirement that 
$[x_{i(t)},x_{j(t)}]$ be a wire of $\Phi $.  
A {\it graph} of $\Phi$ 
(cf.~Fig.~\ref{4FIG7}(2)) is any 1-dimensional polyhedron 
$\gr{\Phi} \sub (\R^2\times\{0\})\cup(\R\times
Y(\Phi)\times[0,\infty[)\sub\R^3$ 
such that 
\begin{romanlist2}{(ii)}
\item $\gr(\Phi)\cap(\R^2\times\{0\})$ is the union of the 
posts of $\Phi\sub\R^2=\R^2\times\{0\}$, 
\item the restriction to the closure in $\R^3$ 
of $\gr(\Phi)\cap(\R\times Y(\Phi)\times ]0,\infty[)$ 
of the projection pr$_{1,2}$ is a homeomorphism onto the union of the 
wires of $\Phi$.
\end{romanlist2}
Any two graphs of $\Phi$ are vertically isotopic with $\Phi$ fixed. 

A {\it charge} on $\Phi$ is a function $\eps:Y(\Phi)\to\{1,-1\}$; 
$(\Phi,\eps )$ is a {\it charged fence}. 
 
\subsection{Constructions with Fences}

\begin{construction} 
\label{band representation to fence}
Given an embedded band representation $\mathbf b$ 
of length $k$ in $B_n$, with $b(t)=:\s_{i(t),j(t)}^{\eps(t)}$ 
for $t=1,\dots,k$, construct a charged fence 
$(\Phi[\mathbf b],\eps[\mathbf b])$ as follows: 
the posts of $\Phi[\mathbf b]$ are $\{s\}\times[1,k]$
for  $s=1,\dots,n$, its wires are $[i(t),j(t)]\times\{t\}$ for $t=1,\dots,k$, 
and the charge is $\eps[\mathbf b](t):=\eps(t)$.  

In particular, if $\mathbf b$ is quasipositive, then $\eps[\mathbf b](t)=1$ 
for all $t$; we will write $\mathbf +$ (rather than 1) for 
this charge.
\end{construction} 

\begin{construction} 
\label{fence to band representation}
Conversely, given a charged fence 
$(\Phi,\eps)$, construct an embedded band representation 
$\mathbf b[\Phi,\eps]$ of length $k:=\card(Y(\Phi))$ 
in $B_{\card(X(\Phi))}$, by setting 
$b[\Phi,\eps](t):=\s_{i(t),j(t)}^{\eps({y_t})}$.
\end{construction} 
 
Constructions~\ref{band representation to fence} and 
\ref{fence to band representation} 
show that embedded band representations 
and charged fences are essentially 
the same: $\mathbf b\mapsto(\Phi[\mathbf b],\eps[\mathbf b])$ and 
$(\Phi,\eps )\mapsto\mathbf b[\Phi,\eps ]$ are mutual inverses 
up to an obvious equivalence relation on charged fences.
 
A fence $\Phi$ is {\it annular} if $\gr(\Phi)$ is a link.
Embellished with over-crossings as in Fig.~\ref{4FIG7}(2), an annular 
fence becomes an unoriented link diagram for the link 
which is its graph; call such a diagram 
an {\it annular fence diagram}.
 
\begin{lemma} (\cite{2,16})  All unoriented links 
have annular fence diagrams. 
\end{lemma}

\begin{construction}
\label{fence to positive plat} Given $\Phi\sub\R^2$,  
an annular fence with respect to coordinates $(x,y)$,  
put $\xi:=x-y$, $\eta:=x+y$, and consider 
the restriction $\eta\vert\Phi$; this has an equal 
number, say $m(\Phi)$, of local maxima (at ``upper right'' corners 
of $\Phi$) and local minima (at ``lower left'' corners).  At 
each local extremum, apply one of the procedures illustrated in 
Fig.~\ref{4FIG8}.
\begin{figure}[htbp]
\setlength{\unitlength}{1in}
\centering
\begin{picture}(4.0,1.5)(0,0)
\includegraphics[width=4.0in]{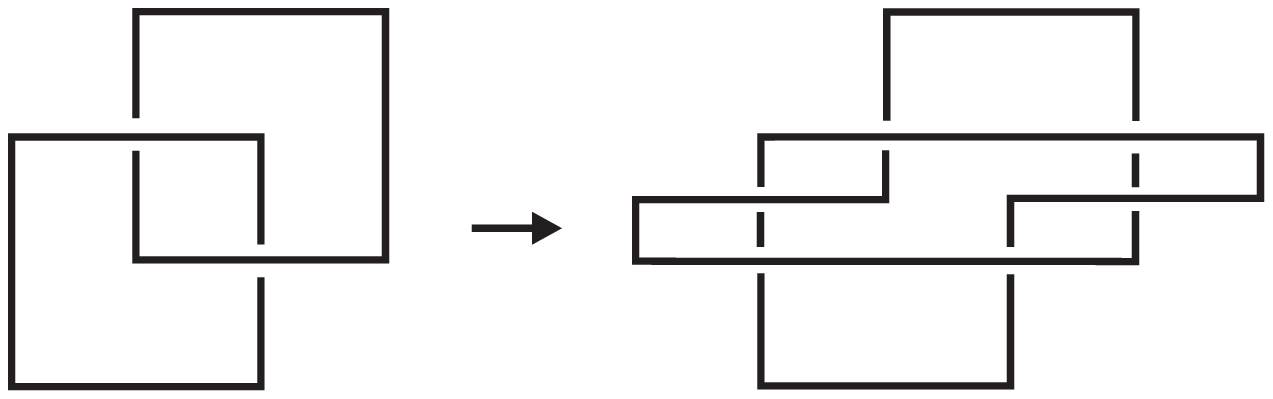}
\end{picture}
\vspace*{13pt}
\fcaption{\label{4FIG8}}
\end{figure}
\noindent 
This replaces $\Phi$ by an annular fence $\Phi'$ 
such that 
\begin{romanlist2}{(iii)}
\item $\gr(\Phi')$ is isotopic to $\gr(\Phi)$ 
(the link diagrams obtained from $\Phi$ and $\Phi'$ 
are regularly homotopic), 
\item $\eta\vert\Phi'$ also has $m(\Phi)$ local maxima 
(and $m(\Phi)$ local minima), and 
\item $\eta\vert\Phi'$ has a single 
local maximum value and a single local minimum value. 
\end{romanlist2}
When we view the link diagram derived from $\Phi'$ 
in $(\xi,\eta)$-coordinates,  
we see the diagram of a positive braidword in $B_{2m(\Phi)}$, say 
$\mathbf p[\Phi]$,  with loose ends joined at the local extrema of 
$Y\vert\Phi'$ according to a suitable plat-plan, say  
$\Pi[\Phi]$, so that $\beta^{\Pi[\Phi]}(\mathbf p[\Phi])$ and 
$\gr(\Phi)$ are ambient isotopic, cf.~Fig.~\ref{4FIG9}.
\end{construction} 
\begin{figure}[htbp]
\setlength{\unitlength}{1in}
\centering
\begin{picture}(2.2,2.0)(-.5,0)
\includegraphics[width=1.5in]{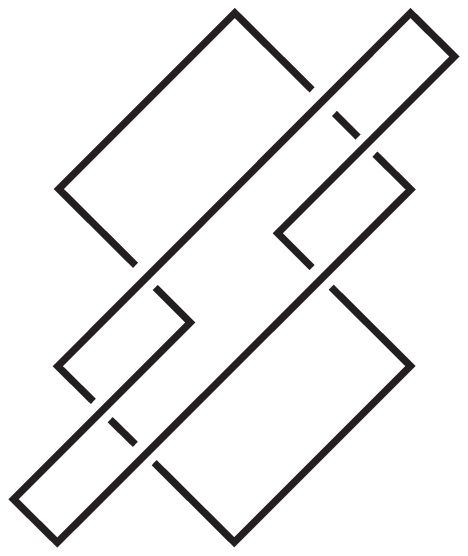}
\end{picture}
\vspace*{13pt}
\fcaption{$\beta^\Pi(\mathbf p)$ for 
$\mathbf p=(\s_2,\s_1,\s_1,\s_3,\s_3,\s_2)$
and $\pi_\sqcup=\pi_\sqcap=(1~2)(3~4)$.\label{4FIG9}}
\end{figure}

\begin{construction}
\label{positive plat to fence} Conversely, given a positive 
braidword $\mathbf p$ in $B_{2m}$ and a plat-plan $\Pi$,  
construct a diagram for $\beta^\Pi(\mathbf p)$ 
such that each segment has slope 
$1$ or $-1$, cf.~Fig.~\ref{4FIG10}; viewed 
in $(\xi,\eta)$-coordinates, the underlying 
graph of this diagram is an annular fence $\Phi[\mathbf p,\Pi]$.
\begin{figure}[htbp]
\vspace*{13pt}
\setlength{\unitlength}{1in}
\centering
\begin{picture}(2.4,2.2)(-.4,0)
\includegraphics[width=2.0in]{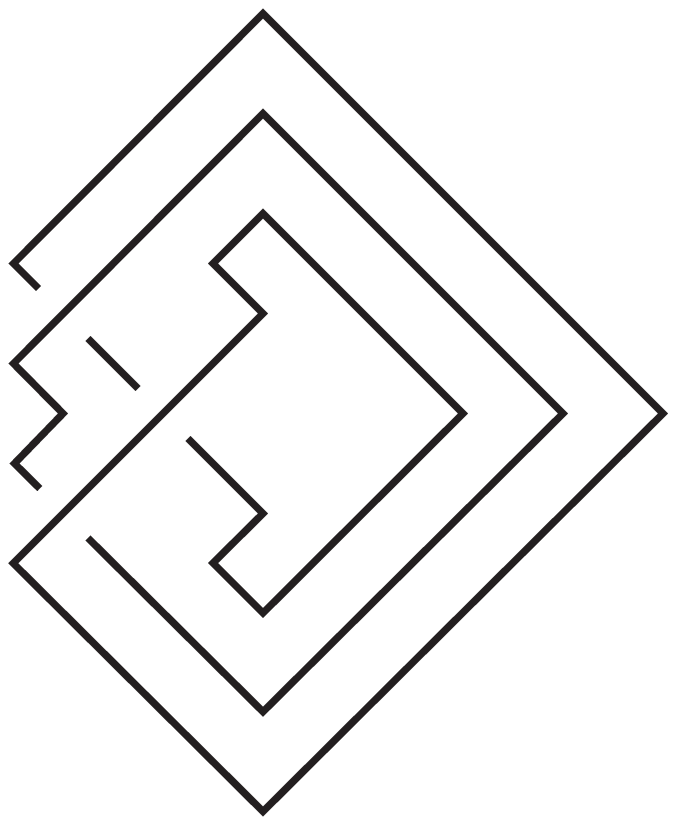}
\end{picture}
\fcaption{\label{4FIG10}}
\end{figure}
\end{construction}

Constructions~\ref{fence to positive plat} 
and \ref{positive plat to fence}
show that annular fence diagrams and 
positive plats are essentially the same: 
$(\mathbf p,\Pi)\mapsto\Phi[\mathbf p,\Pi]$ and 
$\Phi\mapsto(\mathbf p[\Phi],\Pi[\Phi])$ are mutual inverses up 
to obvious equivalences (for annular fence diagrams, the equivalence 
is an adaptation of regular homotopy). 
 
\subsection{Fences and the Modulus of Quasipositivity }

We have seen that an annular embedded band 
representation $\mathbf b$, an annular braided surface $S(\mathbf b)$,  
a charged annular fence, and an appropriately framed positive plat 
all convey the very same information.  Specialize to the case 
that $\mathbf b$ is quasipositive, i.e., $\eps=\mathbf{+}$. 
Nothing is lost by assuming that $S(\mathbf b)$ is a single annulus, 
say $A(K,f)$.  Let $\mathbf p:=\mathbf p[\Phi[\mathbf b]]$,  
$\Pi:=\Pi[\Phi[\mathbf b]]$, and $m:=m(\Phi[\mathbf b])$,  
so that $\beta^\Pi(\mathbf p)$ is a positive plat on $2m$ strings 
that realizes (the unoriented knot underlying) $K$. 
 
\begin{proposition}  In this case, the framing $f$ is equal 
to the writhe of the annular fence diagram of $\Phi[\mathbf b]$ 
(equivalently, the positive plat diagram of $\beta^\Pi(\mathbf p)$)
diminished by $m$.
\end{proposition} 

\begin{figure}[htbp]
\vspace*{13pt}
\setlength{\unitlength}{1in}
\centering
\begin{picture}(1.25,1.25)(0,0)
\includegraphics[width=1.25in]{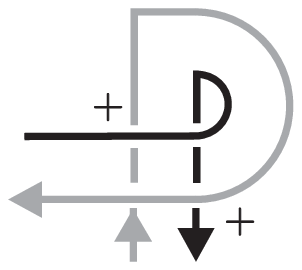}
\end{picture}
\vspace*{13pt}
\fcaption{\label{4FIG11}}
\end{figure}

\begin{proof} It is clear that each crossing in the annular fence
diagram makes the same contribution to $f$ as to the writhe;
Fig.~\ref{4FIG11} shows how each of the $m$ ``upper right'' corners
adds $1$ to the linking of two components of
the boundary of the annulus, and thus diminishes the framing by 1. 
\end{proof}  

\begin{corollary}
\label{q from positive plats}
$q(K)$ is the maximum, over 
all $m$ and all realizations of $K$ as a plat $\beta^\Pi(\mathbf p)$ 
of a positive braidword $\mathbf p$ in $B_{2m}$, of the writhe of 
$\beta^\Pi(\mathbf p)$ diminished by $m$.  
\end{corollary}

\begin{theorem}
\label{bound q from below}
If $\beta^\Pi(\mathbf p)$ 
is any realization of $K$ as a positive $\Pi$-plat on $2m$ strings,
then $q(K)$ is greater than or equal to the writhe of 
$\beta^\Pi(\mathbf p)$ diminished by $m$; equivalently, if $\gr(\Phi)$ 
is any realization of $K$ as the graph of an annular fence, then 
$q(K)$ is greater than or equal to the writhe of the annular fence 
diagram of $\Phi$ diminished by $m(\Phi )$. 
\end{theorem}

\begin{corollary}
\label{q of closed positive braid--inequality} 
(\cite{3})  If $\mathbf p$ is a 
positive braidword in $B_n$,  
then $q(\widehat{\beta}(\mathbf p))\geq e(\beta(\mathbf p))-n$. 
\end{corollary}
 
\begin{proof} As is well-known, if $\beta\mapsto\beta^{(n)}$ denotes 
the injection $B_n\to B_{2n}$ which takes $\s_i\in B_n$ to 
$\s_i\in B_{2n}$, and 
$\pi_\sqcup=\pi_\sqcap=(1\ 2n)(2\ 2n-1)\dotsm(n\ n+1)$,  
then for all $\beta\in B_n$, the closed braid $\widehat{\beta}$ 
and the $\Pi$-plat $(\beta^{(n)})^\Pi$ are ambient isotopic. 
\end{proof} 
 
\begin{proposition}\label{subadditivity of q+1} 
For any knots $K_k$, 
$q(\Connsum_k K_k)+1\geq \sum _k(q(K_k)+1)$.   
\end{proposition}

\begin{proof} If $\Phi$ is any annular fence, then the top of the rightmost 
post of $\Phi$ is a local maximum of $\eta\vert\Phi$, and no 
point of the leftmost post of $\Phi$ is a local maximum of 
$\eta\vert\Phi$.  Let $\Phi_k$ be an annular fence such that 
$gr(\Phi_k)=K_k$ and the writhe of the annular fence diagram of $\Phi_k$ 
is $q(K_k)+m(\Phi_k)$.  By applying appropriate translations and 
homotheties, we may assume that the rightmost post of $\Phi_k$ 
is the leftmost post of $\Phi_{k+1}$ for $k=1,\dots,N-1$.  The fence 
obtained from $\bigcup_k\Phi_k$ by deleting all the interiors of 
the common posts (which may be suggestively denoted by $\Connsum_k\Phi_k$)
is annular, the writhe of its fence diagram is the sum of the writhes of 
the fence diagrams of the $\Phi_k$, 
$m(\Connsum_k\Phi_k)-1=\sum_k (m(\Phi_k)-1)$, and
$\gr(\Connsum_k\Phi_k)=\Connsum_k\gr(\Phi_k)$; so the proposition 
follows from Theorem~\ref{bound q from below}.
\end{proof} 
 
\begin{remark}\label{strictness?} 
I do not know if the inequality in 
Proposition~\ref{subadditivity of q+1} is ever strict.  
\end{remark} 

\noindent
{\bf Historical remark.} Fences are my synthesis of (i) some diagrams 
that H. Morton used to describe certain Hopf-plumbed fiber surfaces 
in 1982 at Les-Plans-sur-Bex, Switzerland, and 
(ii) ``square bridge projections'' as described by H. Lyon in 1977 
in Blacksburg, Virginia (\cite{16}): an unoriented link is 
in ``square bridge position'' if and only if it is the graph 
of an annular fence (see below).  Square bridge 
projections have frequently been rediscovered---for instance by 
Thurston, Erlandsson \cite{17}, and Kuhn \cite{18}, who (jointly 
and severally) call them ``barber-pole projections''.  

\section{Upper Bounds for the Modulus of Quasipositivity}

In this section we derive various upper bounds for $q(K)$ 
(Corollaries~\ref{R-bound on q}, 
\ref{framed-polynomial bound on q}, 
\ref{F-polynomial bound on q}, 
and \ref{G-polynomial bound on q}), 
all given in terms of link polynomials
and based on a fundamental result of 
Morton and Franks~\&~Williams.
 
\begin{theorem}\label{MFW theorem} (\cite{14}, \cite{19})
For all $n$, for all 
$\beta\in B_n$,  we have $\ord_v P_{\widehat\beta}\geq e(\beta)-n+1$.
\end{theorem}

\begin{corollary}\label{v-order of qp annular boundary}
 (\cite{9}, \cite{20})  If $\mathbf b$ 
is a quasipositive annular embedded band representation, then 
$\ord_v P_{\widehat\beta(\mathbf b)}\geq 1$.  Equivalently, 
if $A(L,f)$ is quasipositive, then 
$\ord_v P_{\Bd A\{L,f\}}\geq 1$.  
\end{corollary}
 
\begin{corollary}\label{R-bound on q} 
For any knot $K$, $q(K)\leq -1+\ord_v R_K$.  
\end{corollary}
 
\begin{proof}  For all $L$, $\ord_v P_L\leq\ord_v R_L$.  By
Corollary~\ref{R of disjoint sum}, 
$R_{\Bd A(K,f)}\negthinspace=\negthinspace%
(v^{-1}-v)v^{-2f}R_K(v)^2$ for all $(K,f)$.
Let $\Bd A(K,f)$ be quasipositive; then 
$1\leq\ord_v P_{\Bd A(K,f)}\leq
\ord_v  
R_{\Bd A(K,f)}$ %
$\leq -1-2f+2\ord_v R_K \leq -1+2\ord_v R_K-2q(K)$ 
by Corollary~\ref{v-order of qp annular boundary}. 
\end{proof} 
 
\begin{theorem} If $A(L,f)$ is quasipositive, then 
$\ord_v(\{L,f\}-(-1)^{c(L)})\geq0$. 
\end{theorem}

\begin{proof} When $A(L,f)$ is quasipositive, so is
$A(L',f\vert L')$ for each non-empty sublink 
$L'\sub L$; then
$$
\ord_v (\{L,f\}-(-1)^{c(L)})
\geq\min\limits_{\emptyset\ne L'\sub L} 
\ord_v\left((v^{-1}-v)z^{-1}P_{\Bd A({L'},
f\vert{L'})}\right)\geq 0
$$
by Corollary~\ref{v-order of qp annular boundary}.
\end{proof}  
	
\begin{corollary}
\label{framed-polynomial bound on q} 
(\cite{9,20}) For any knot $K$, 
$q(K)\leq\frac{1}{2}\ord_v\{K,0\}$.  
\end{corollary}

\begin{theorem} If $A(L,f)$ is quasipositive, then 
$\deg_a F^*_{L}\leq-1-\phi(L,f)-2\tau(L)$. 
\end{theorem} 

\begin{proof} If $A(L,f)$ is quasipositive, then 
\begin{align*}
0
&\le\ord_v(\{L,f\}-(-1)^{c(L)})\\
&\le\ord_v\big((1+\{L,f\})
\negmedspace\negmedspace\negmedspace
 \pmod 2\big)\\
&\le\ord_v(1+v^{-4\tau(L)-2\phi(L,f)}(1+(v^{-2}+v^2)z^2)F^*_L(v^{-2},z^2)).
\end{align*}
But if $\deg_a(F^*_{L})>-1-\phi(L,f)-2\tau(L)$, then
\begin{align*}
0
&>-2-2\phi(L,f)-4\tau(L)-2\deg_a(F^*_{L})\\
&=\ord_v(v^{-4\tau(L)-2\phi(L,f)}(1+(v^{-2}+v^2)z^2)F^*_L(v^{-2},z^2))\\
&=\ord_v(1+v^{-4\tau(L)-2\phi(L,f)}(1+(v^{-2}+v^2)z^2)F^*_L(v^{-2},z^2))
\end{align*}
since ${\ord_v 1}=0$ and $\ord_v(A+B)=\min\{\ord_v A, \ord_v B\}$ 
if $\ord_v A \ne\ord_v B$.\end{proof}
 
\begin{corollary}\label{F-polynomial bound on q} 
(\cite{9,20}) For any knot $K$, 
$q(K)\leq-1-\deg_aF^*_K$.  
\end{corollary}
 
\begin{corollary}\label{G-polynomial bound on q} 
For any knot $K$,  for 
$k=0,1$,  $q(K)\leq-1-\deg_a G^k_K$.  
\end{corollary}
 
\begin{proof} This follows from 
Corollary~\ref{F-polynomial bound on q} since 
$\deg_a F^*_{L}=\max\{\deg_a G^0_L, \deg_a G^1_L\}$.  (The 
case $k=0$ also follows from Theorem~\ref{congruence theorem}
and Corollary~\ref{R-bound on q}.) 
\end{proof} 
 
\section{Some Calculations of the Modulus of Quasipositivity}

\subsection{Summary of results}

To illustrate the usefulness of the upper and lower bounds derived
in the preceding sections, we compute $q(K)$ exactly for various 
infinite classes of knots.  (The author and Michel Boileau have made use 
of these calculations in connection with constructions of Stein manifolds
with pseudoconvex boundaries of interesting topological types.)

\subsection{Positive Closed Braids}

\begin{lemma}\label{positive braidwords}(\cite{1})
Let $\mathbf p$ be a positive 
braidword in $B_n$.  Either no generator $\s_i$ of $B_n$ appears 
more than once in $\mathbf p$,  or there is a positive braidword $\mathbf q$ 
in $B_n$ such that $q(1)=q(2)$, $e(\beta(\mathbf q))=e(\beta(\mathbf p))$, 
and $\widehat{\beta}(\mathbf q)$ is isotopic to 
$\widehat{\beta}(\mathbf p)$.
\end{lemma}

\begin{theorem}
\label{R of positive braid}
If $\mathbf p$ is a positive braidword in $B_n$,  
then $\ord_v R_{\widehat\beta(\mathbf p)}=e(\beta(\mathbf p))-n+1$, 
and the coefficient of $v^{e(\beta(\mathbf p))-n+1}$ in 
$R_{\widehat\beta(\mathbf p)}$ is a positive integer.
\end{theorem}
 
\begin{proof}  If $n=1$ then $\mathbf p$ is empty, $\widehat{\beta}(\mathbf p)$ 
is an unknot, and the conclusion holds for $\mathbf p$.  
Let $n>1$, and assume the conclusion for all positive braidwords on fewer than 
$n$ strings, and all positive braidwords on $n$ strings with exponent 
sum less than $e(\beta(\mathbf p))$.  If $\widehat\beta(\mathbf p)$ 
has more than one component, let 
$\emptyset\ne L'\sub\widehat{\beta}(\mathbf p)$,
$\emptyset\ne L'':=\widehat\beta(\mathbf p)\setminus L'$; then there 
are positive braidwords $\mathbf p'$, $\mathbf p''$ in $B_n'$ and $B_n''$
(where $n'+n''=n$), with $L'=\widehat\beta(\mathbf p')$,
$L''=\widehat\beta(\mathbf p'')$, 
$2\link(L',L'')=e(\beta(\mathbf p)-e(\beta(\mathbf p')-e(\beta(\mathbf p'')$;
by the inductive hypothesis and 
Corollary~\ref{R of disjoint sum}, the conclusion 
holds for $\mathbf p$.  If $\widehat\beta(\mathbf p)$ 
has exactly one component,
then one of the alternatives in Lemma~\ref{positive braidwords} is the case.
If no generator appears more than once, then each generator 
appears exactly once, so $\widehat\beta(\mathbf p)=O$ and the conclusion holds 
for $\mathbf p$.  If a generator appears more 
than once, then let $\mathbf q$ be as in 
Lemma~\ref{positive braidwords}, $\mathbf q' :=
(q(2),\dots,q(e(\beta(\mathbf p))))$, 
$\mathbf q'':=(q(3),\dots,q(e(\beta(\mathbf p))))$; 
by Lemma~\ref{R-polynomial}, $R_{\widehat\beta(\mathbf p)}= 
R_{\widehat\beta(\mathbf q)}= 
vR_{\widehat\beta(\mathbf q')}
+v^2 R_{\widehat\beta(\mathbf q'')}$ so by the inductive hypothesis 
the conclusion holds for $\mathbf p$.  
\end{proof}  
 
\begin{corollary}
\label{q of closed positive braid}
Let $\mathbf p$ be a positive 
braidword in $B_n$ such that $\widehat\beta(\mathbf p)$ 
is a knot. Then $q(\widehat\beta(\mathbf p))=e(\beta(\mathbf p))-n$.  
\end{corollary} 

\begin{proof} Combine Theorem~\ref{R of positive braid}
with Corollaries 6 and 8. 
\end{proof} 
 
\begin{remark} In terms of the Milnor number of the 
fibered knot $\widehat\beta(\mathbf p)$, 
Corollary~\ref{q of closed positive braid} says that
$q(\widehat\beta(\mathbf p))=\mu(\widehat\beta(\mathbf p))-1$.  The appearance 
of $\mu$ is undoubtedly a red herring; there may be some generalization 
of Corollary~\ref{q of closed positive braid}
to other fibered knots, but my guess is that it is at best 
an inequality, and involves not only $\mu$ but also the {\it enhancement}
$\lambda$ (which happens to equal $0$ for a positive 
closed braid: $\lambda$ is sensitive to handedness, as $q$ seems 
to be and $\mu$ manifestly is not), cf.~\cite{20} and references cited 
therein.
\end{remark} 
 
\subsection{Two-strand Torus Knots}

\begin{theorem}  $q(O\{ 2,2k+1\})$ equals $2k-1$ for $k\geq 0$, $-1$ 
for $k=-1$, and $4k+2$ for $k\leq-2$.
\end{theorem}

\begin{proof} The cases with $k\geq-1$ are already done: if $k\geq 0$, 
then represent $O\{ 2,2k+1\}$ as the closed positive $2$-string braid 
$\widehat{\s_1^{2k+1}}$ and apply 
Corollary~\ref{q of closed positive braid}; if $k=-1$, 
then note that $O\{2,-1\}=O$.  

Let $k\leq -2$.  Put $r:={-(2k+1)}\geq 3$ and represent 
$O\{2,2k+1\}$ as the positive $\Pi$-plat 
$(\s_1\s_3\dots\s_r)^\Pi$ on $2r$ strings, 
where $\pi_\sqcup=\pi_\sqcap=(1~2r)(2~3)\dotsm({2r-2}~{2r-1})$.
Then the writhe of the plat diagram is $-r$, so 
by Theorem~\ref{bound q from below},
$q(O\{2,2k+1\})\geq -r-\frac{1}{2}(2r)=4k+2$.  
To prove the opposite inequality 
and finish the proof, it suffices 
by Corollary~\ref{G-polynomial bound on q}
to show that $\deg_a G^1_{O\{2,2k+1\}}\geq-4k-3$ for $2k+1\leq -3$.  
By Lemma~\ref{G-polynomial}, $G^1_{O\{2,-3\}}=a^2+a^3+a^5$,
$G^1_{O\{2,-4\}}=a^3+a^6+a^7$, 
$G^1_{O\{2,-5\}}=a^5+a^8+a^9$,  and 
$G^1_{O\{2,m\}}=a^3G^1_{O\{2,m+3\}}+a^{-2m-3}+a^{-2m-1}$ for 
$m\leq -6$; by induction we have $\deg_a G^1_{O\{2,m\}}=-2m-1$ 
for $m\leq-3$.  
\end{proof} 
 
\subsection{Some Positive Pretzel Knots}

Let $\pi_\sqcup=\pi_\sqcap=(1~6)(2~3)(4~5)\in\Ss_6$.  
The unoriented link $\Ps(r,s,t):=(\s^r_1\s^s_3\s^t_5)^\Pi$ 
is called a {\it pretzel}; it is {\it positive} iff $r,s$, and $t$ are all 
non-negative.  Note that $\Ps(r,s,t)$ is a knot iff two or 
three of $r,s,t$ are odd, and that $\Ps(r,s,t)$ is ambient isotopic 
to $\Ps(s,t,r)$. 
 
\begin{theorem} If $r,s,t+1\geq 1$ are odd, then $q(\Ps(r,s,t))=-3+r+s-t$. 
\end{theorem} 

\begin{proof} In the braidword diagram of 
$(\s_1,\dots,\s_1,\s_3,\dots,\s_3,\s_5,\dots,\s_5)$, 
in this case, 
the crossings that correspond to $\s_1$ 
and $\s_3$ 
are positive and those that correspond to 
$\s_5$ are negative, so the writhe of the associated knot diagram 
of $\Ps(r,s,t)$ is $r+s-t$; by 
Theorem~\ref{bound q from below}, $q(\Ps(r,s,t))\geq-3+r+s-t$.  

To establish the opposite inequality, it suffices, 
by Corollary~\ref{R-bound on q}, to show that 
$\ord_v R_{\Ps(r,s,t)}=-2+r+s-t$.  If $t=0$, then 
$\Ps(r,s,0)=O\{2,r\}\connsum O\{2,s\}$, while if $t\geq 2$, then 
(by considering any one of the negative crossings) $R_{\Ps(r,s,t)}=
v^{-2}R_{\Ps(r,s,t-2)}-v^{-1}R_{O\{2,r+s)}$; easy inductions complete the proof. 
\end{proof} 
 
\begin{remark} A similar calculation shows that 
$-3-r-s-t\leq q(\Ps(r,s,t))\leq -2-r-s-t$ when $r,s,t\geq 1$ are all odd.
\end{remark} 

\baselineskip 13pt minus1pt
\section{The Modulus of Quasipositivity of a Slice Knot}

Kronheimer~\&~Mrowka \cite{22} have recently announced a gauge-theoretic 
proof of the following long-conjectured result.

\begin{theorem}
\label{KM theorem}
Let $\Gamma\sub\C^2$ be a smooth complex-algebraic curve.
If $\Gamma$ intersects 
$S^3:=\{(z,w)\in\C^2:\left|z\right|^2+\left|w\right|^2=1\}$
transversely, then no smooth orientable surface 
$S\sub D^4:=\{(z,w)\in\C^2:\left|z\right|^2+\left|w\right|^2\le1\}$
without closed components,
such that $\Bd S= \Gamma\cap S^3$, has larger Euler characteristic 
than $\Gamma\cap D^4$.
\end{theorem}

A knot $K\sub S^3$ is {\it slice} if $K=\Bd D$ for some smooth
$2$-disk $D\sub D^4$.  

\begin{proposition}
\label{q(K)<0 for K slice} 
If $K$ is a slice knot, then $q(K)<0$.
\end{proposition}

\begin{proof}  Let $K$ be a knot.  If $K$ is slice,
then $\Bd A(K,0)$ bounds a surface in $D^4$ of Euler characteristic two
(namely, the union of two disjoint smooth $2$-disks).  
On the other hand, if $q(K)\ge 0$, then 
(by Proposition~\ref{quasipositive annuli}) 
the annulus $A(K,0)$ is quasipositive, and it follows from \cite{1} 
(cf.~also \cite{21}) that there exists 
a smooth complex-algebraic curve $\Gamma$ such that 
$\Gamma\cap S^3$ is a link of type $\Bd A(K,0)$ (the intersection being
transverse), while $\Gamma\cap D^4$ is a surface of Euler characteristic zero 
(namely, a ``push-in'' of $A(K,0)$).  Now the proposition follows from 
Theorem~\ref{KM theorem}.
\end{proof}

\begin{remark} Another interesting consequence of 
Theorem~\ref{KM theorem} is that,
for quasipositive knots, slice implies ribbon.  
\end{remark} 

\nonumsection{Acknowledgments}
This research was partially supported by NSF grant DMS-8801959.

\nonumsection{Addendum (December 2001)}

Bennequin's proof \cite{23} that the {\it maximal 
Thurston--Bennequin invariant $\TB(K)$} of a knot in $S^3$ 
is an integer (rather than $\infty$) sparked considerable
interest in finding ways to compute $\TB(K)$, or at least
bound it above (see, e.g., \cite{24}).  

In \cite{25}, it was shown that $\TB(K)$
is identical to the modulus of quasipositivity $q(K)$; thus
the various bounds on, and calculations of, $q$ derived above
are equally bounds on, or calculations of, $\TB$.  Similar 
(sometimes sharper) results for $\TB$ have been derived using 
a variety of different methods by a number of researchers, 
among them Fuchs \& Tabachnikov \cite{26}, 
Tabachnikov \cite{27},
Epstein \cite{28},
Chmutov and Goryunov \cite{29},
Kanda \cite{30},
Tanaka \cite{31},
Goryunov and Hill \cite{32},
Etnyre and Honda \cite{33},
Ng \cite{34},
and Ferrand \cite{35}. 

Torisu \cite{36} and Etnyre \& Honda \cite{37} 
have announced that $\TB+1$ is additive; in light of \cite{25}, 
this settles the issue raised in Remark~\ref{strictness?}, 
by showing that $q+1$ is additive (the inequality in 
Proposition~\ref{subadditivity of q+1} may be replaced
by an equation).

\nonumsection{References}

\end{document}